\documentclass[12pt,notitlepage]{amsart}
\usepackage{latexsym,amsfonts,amssymb,amsmath,amsthm}
\usepackage{color}
\newcommand{\mz}{\ensuremath{\mathbb Z}}
\newcommand{\mr}{\ensuremath{\mathbb R}}

\newcommand{\mc}{\ensuremath{\mathbb C}}

\newcommand{\half}{\ensuremath{ \frac{1}{2}}}

\newcommand{\intR}{\int_{-\infty}^{\infty}}

\newcommand{\leg}[2]{\left(\frac{#1}{#2}\right)}
\newcommand{\e}[2]{e\left(\frac{#1}{#2}\right)}

\theoremstyle{plain}		
	\newtheorem{mytheo}{Theorem}[section]

	\newtheorem{mycoro}[mytheo]{Corollary}
     \newtheorem{mylemma}[mytheo]{Lemma}

\theoremstyle{remark}

\numberwithin{equation}{section}
\begin{document}
\author{Xiaoqing Li}
\address{Department of Mathematics \\
State University of New York at Buffalo \\
Buffalo, NY, 14260
}
\email{XL29@buffalo.edu}

\author{Matthew P. Young} 
\address{Department of Mathematics \\
	  Texas A\&M University \\
	  College Station \\
	  TX 77843-3368}
\email{myoung@math.tamu.edu}
\thanks{This material is based upon work supported by the National Science Foundation under agreement Nos. DMS-0901035 (X.L.), DMS-0758235 (M.Y.),  and DMS-0635607 (X.L. and M.Y.).  Any opinions, findings and conclusions or recommendations expressed in this material are those of the authors and do not necessarily reflect the views of the National Science Foundation.}

% \begin{abstract}
% \end{abstract}
%\today
\title{Additive twists of Fourier coefficients of symmetric-square lifts}
\maketitle

\section{Introduction}
It is a classical problem to estimate sums involving the Fourier coefficients of a modular form.  For instance, if $f(z) = \sum_{n} \lambda(n) n^{\frac{k-1}{2}} e(nz)$ is a weight $k$ holomorphic cusp form then $\sum_{n \leq N} \lambda(n) e(\alpha n) \ll_f N^{1/2} \log{2N}$, uniformly in $\alpha \in \mr$; see Theorem 5.3 of \cite{IwaniecTopics}.  This amounts to ``square-root'' cancellation, uniformly in $\alpha$.  This uniformity in $\alpha$ is pleasant in applications and allows one to study sums over $a(n)$ with $n$ lying in a fixed arithmetic progression; e.g., see Corollary 5.4 of \cite{IwaniecTopics}.

In a different direction, it is interesting to study sums of Fourier coefficients over other sequences, such as the values of a quadratic polynomial \cite{Blomer}.  In particular, sums involving $\lambda(n^2)$ appear in many applications since these appear in the Dirichlet series coefficients of the symmetric-square $L$-function.  In this paper, we study the following sum
\begin{equation}
 \sum_{n \leq N} \lambda_j(n^2) e(\alpha n),
\end{equation}
where $\lambda_j(m)$ is the $m$-th Hecke eigenvalue of a Maass form for the full modular group with Laplacian eigenvalue $1/4 + t_j^2$.  It is well-known that the symmetric-square lift of a $GL_2$ Maass form is a $GL_3$ Maass form \cite{GJ}.  Therefore, the above sum is closely related to the following
\begin{equation}
\label{eq:basicsum}
S_F(N)= \sum_{n \leq N} A_F(1,n) e(\alpha n),
\end{equation}
where $A_F(1,n) = \sum_{ml^2 = n} \lambda_j(m^2)$ and
\begin{equation}
\label{eq:Ldef}
 L(F,s) = \sum_{n=1}^{\infty} A_F(1,n) n^{-s}
\end{equation}
is the $L$-function associated to a $SL_3(\mz)$ Maass form.  S.D. Miller \cite{Miller} showed that $S_F(N) \ll_F N^{3/4+\varepsilon}$
uniformly in $\alpha \in \mr$ but with an implied constant depending on $F$.  A key tool is the $GL_3$ Voronoi formula proved by Miller-Schmid \cite{MS}.

In this note we consider the dependence of \eqref{eq:basicsum} on the form $F$, that is, in terms of the spectral parameter $t_j$.  Our motivation arises from some recent studies in the analytic theory of automorphic forms that require estimates which are uniform with respect to the automorphic form.  For example, \cite{LY} required an application of the Voronoi formula with a varying underlying form.  In the more well-known $GL_2$ case we have available uniform asymptotic expansions of the Bessel functions which unfortunately are not known for the Whittaker functions in the $GL_3$ case.  The problem of estimating \eqref{eq:basicsum} with a varying form $F$ is attractive for a few reasons.  For one, it is fundamental to understand correlations of the $GL_3$ Fourier coefficients with a linear phase.  
Secondly, the limited number of parameters ($F$, $N$, $\lambda$) apparently makes a good setting for exploring new behavior in these types of sums.  Finally,  it turns out that there is some unexpected new behavior not present in the case of $F$ fixed.

The analysis of \eqref{eq:basicsum} leads to some difficult technical problems in the theory of exponential integrals.  Restricting our considerations to symmetric-square lifts serves as a relatively nice compromise between generality and difficulty; already this case is quite tricky to analyze.  It would be interesting to study \eqref{eq:basicsum} for a general Maass form $F$.

Our main result is
\begin{mytheo}
\label{thm:Millertype}
Suppose $F$ is the symmetric-square lift of a $SL_2(\mz)$ Hecke-Maass form with $A_F(1,1) = 1$.  Then assuming the Ramanujan conjecture at the finite places for $F$, we have with $D=1/4$,
\begin{equation}
\label{eq:mainthmRama}
 \sum_{n \leq N} A_F(1,n) e(\alpha n) \ll N^{3/4 + \varepsilon} \lambda_F(\Delta)^{D + \varepsilon},
\end{equation}
 where $\lambda_F(\Delta)$ is the Laplace eigenvalue of $F$.  The implied constant depends on $\varepsilon > 0$ only. 
Unconditionally, we have \eqref{eq:mainthmRama} with $D=1/3$. 
% \begin{equation}
%  \label{eq:mainthmNoRama}
%  \sum_{n \leq N} A_F(1,n) e(\alpha n) \ll N^{3/4 + \varepsilon} \lambda_F(\Delta)^{D + \varepsilon},
% \end{equation}

\end{mytheo}
If $F$ has Langlands parameters $(iT, 0, -iT)$ then $\lambda_F(\Delta) = 1+T^2$.  The analytic conductor of $L(F,1/2)$ is then $\asymp T^2$. If $F$ is the symmetric-square lift of a Maass form with Laplace eigenvalue $1/4 + t_j^2$ then $T = 2t_j$.

Our approach does not need the full strength of the Ramanujan conjecture.  Instead, we require bounds on sums of the form
\begin{equation}
 \sum_{A \leq n \leq A + B} |A_F(1,n)|,
\end{equation}
where $B \leq A$.  It is even tricky to see exactly how small we require for $B$.

We deduce Theorem \ref{thm:Millertype} from the following
\begin{mytheo}
 \label{thm:mainthmsmooth}
 Let conditions be as in Theorem \ref{thm:Millertype}.
Suppose $w$ satisfies
\begin{equation}
\label{eq:wbounds}
\begin{cases} 
\text{$w$ is smooth with support in the dyadic interval $[N, 2N]$} \\
|w^{(j)}(y)| \leq c_j N^{-j},
\end{cases}
\end{equation}
for some positive real numbers $c_j$, and all $j=0,1,2, \dots$.  Then
\begin{equation}
\label{eq:smoothsum}
 \sum_{n=1}^{\infty} A_F(1,n) e(\alpha n) w(n) \ll N^{3/4 + \varepsilon} \lambda_F(\Delta)^{D+\varepsilon},
\end{equation}
where $D=1/4$ if the Ramanujan conjecture holds, and $D=1/3$ unconditionally.
The implied constant depends on the constants $c_j$ in \eqref{eq:wbounds} and on $\varepsilon > 0$ only.  
\end{mytheo}
We sketch this (fairly standard) deduction in Section \ref{section:unsmoothing} below.

\section{Acknowledgements}
The authors are indebted to Soundararajan for a suggestion which improved Theorem \ref{theo:Molteni}.  We also thank S.D. Miller for encouragement.

\section{Prerequisites on $GL_3$ Maass forms}
We very briefly state the necessary notions of $GL_3$ Maass forms.  We work almost exlusively on the level of $L$-functions so we simply refer to Goldfeld's book \cite{Goldfeld} for the automorphic picture.  Suppose that $F$ is a Maass form of type $(\nu_1, \nu_2)$ for $SL_3(\mz)$ which is an eigenfunction of all the Hecke operators with Fourier coefficients $A_F(m,n)$ normalized so that $A_F(1,1) = 1$ and $|A_F(1,p)| \leq 3$ is implied by the Ramanujan conjectures.  Then the dual of $F$, denoted $\widetilde{F}$ is of type $(\nu_2, \nu_1)$ with $(m,n)$th Fourier coefficient equal to $A_F(n,m) = \overline{A_F(m,n)}$.  It is often convenient to work with the Langlands parameters $(\alpha_1, \alpha_2, \alpha_3)$ defined by
\begin{align}
 \alpha_1 &= -\nu_1 - 2 \nu_2 + 1, \\
 \alpha_2 &= -\nu_1 + \nu_2, \\
 \alpha_3 &= 2 \nu_1 + \nu_2 - 1.
\end{align}
Notice that the dual of $F$ has Langlands parameters $(-\alpha_3, -\alpha_2,-\alpha_1)$.  Furthermore, observe that $\alpha_1 + \alpha_2 + \alpha_3 = 0$.
The $L$-function associated to $F$ defined by \eqref{eq:Ldef} satisfies the functional equation (see Theorem 6.5.15 of \cite{Goldfeld} for example)
\begin{equation}
 G_{\alpha_1, \alpha_2, \alpha_3}(s) L(F,s) = G_{-\alpha_1, -\alpha_2, -\alpha_3}(1-s) L(\widetilde{F}, 1-s),
\end{equation}
where
\begin{equation}
 G_{\alpha_1, \alpha_2, \alpha_3}(s) = \pi^{-3s/2} \Gamma\big(\frac{s-\alpha_1}{2} \big) \Gamma\big(\frac{s-\alpha_2}{2} \big)\Gamma\big(\frac{s-\alpha_3}{2} \big).
\end{equation}

In the special case that $F$ is the symmetric-square lift of a $SL_2(\mz)$ Maass form then the Langlands parameters take the form $\alpha_1 = iT$, $\alpha_2=0$, $\alpha_3 = -iT$.

We quote our key tool, the $GL_3$ Voronoi formula.
Suppose $k=0$ or $1$, and $\psi(x)$ is a smooth, compactly-supported function on the positive reals.  Define
\begin{equation}
\label{eq:psidef}
 \widetilde{\psi}(s) = \int_0^{\infty} \psi(x) x^{s} \frac{dx}{x}.
\end{equation}

For $\sigma > -1 + \max\{-\text{Re}(\alpha_1), -\text{Re}(\alpha_2),  -\text{Re}(\alpha_3)\}$, define
\begin{equation}
\label{eq:psikdef}
 \Psi_k(x) = \frac{1}{2 \pi i} \int_{(\sigma)} (\pi^3 x)^{-s} \frac{\Gamma\left(\frac{1+s+\alpha_1+k}{2}\right)\Gamma\left(\frac{1+s+\alpha_2+k}{2}\right)\Gamma\left(\frac{1+s+\alpha_3+k}{2}\right)}{\Gamma\left(\frac{-s-\alpha_1+k}{2}\right)\Gamma\left(\frac{-s-\alpha_2+k}{2}\right)\Gamma\left(\frac{-s-\alpha_3+k}{2}\right)} \widetilde{\psi}(-s) ds.
\end{equation}
Then define
\begin{gather}
 \Psi_{+}(x) = \frac{1}{2 \pi^{3/2}} ( \Psi_0(x) + \frac{1}{i} \Psi_1(x) ) \\
\label{eq:Psi-}
\Psi_{-}(x) = \frac{1}{2 \pi^{3/2}} (\Psi_0(x) - \frac{1}{i} \Psi_1(x)).
\end{gather}
\begin{mytheo}[\cite{MS}]
\label{thm:Voronoi}
 Let $\psi(x)$ be smooth and compactly-supported on the positive reals.  Suppose $d, \overline{d}, c \in \mz$ with $c \neq 0$, $(c,d) = 1$, and $d \overline{d} \equiv 1 \pmod{c}$.  Then
\begin{multline}
\label{eq:Voronoi}
 \sum_{n > 0} A_F(1,n) \e{n \overline{d}}{c} \psi(n) = c \sum_{n_1 | c} \sum_{n_2 > 0} \frac{A_F(n_2, n_1)}{n_1 n_2} S(d, n_2; c/n_1) \Psi_{+} \leg{n_2 n_1^2}{c^3 }
\\
+ c \sum_{n_1 | c} \sum_{n_2 > 0} \frac{A_F(n_2, n_1)}{n_1 n_2} S(d, -n_2; c/n_1) \Psi_{-} \leg{n_2 n_1^2}{c^3 },
\end{multline}
where $S(a,b;c)$ is the usual Kloosterman sum.
\end{mytheo}

% The analytic conductor of $L(F,1/2)$ (see \cite{IK} for the definition), denoted $C(F)$ here, is $\asymp (1 + |\alpha|)(1+ |\beta|)(1+|\gamma|)$.  It is known that $|\text{Re}(\alpha)| \leq \frac25$ (and similarly for $\beta, \gamma$) [make ref].  Write $\alpha = \alpha_0 + i \alpha_1$ and similarly for $\beta$ and $\gamma$.  The Ramanujan conjecture at the infinite place means $\alpha_0 = \beta_0 = \gamma_0 = 0$.  Suppose that $\alpha_1 \geq \beta_1 \geq \gamma_1$ so that $\alpha_1 > 0 > \gamma_1$ and furthermore
% \begin{equation}
%  -\half \gamma_1 \leq \alpha_1 \leq -2 \gamma_1.
% \end{equation}
% If $|\alpha_1| \asymp |\gamma_1| \asymp T$ then $C(F) \asymp T^2 (1 + |\beta_1|)$.  The local Weyl law implies that ``typically'' $\beta_1 \gg T^{1-\varepsilon}$; in practice this case is relatively easy and in fact for such forms we may remove the Ramanujan conjecture from Theorem \ref{thm:Millertype}.  
% We have (? I don't recall where this came from)
% \begin{mytheo}
% \label{thm:Millertype2}
%  Suppose $F$ is a $GL_3$ Maass form with Langlands parameters $\alpha, \beta, \gamma$.  Then
% \begin{equation}
%  \sum_{n \leq N} A_F(1,n) e(\lambda n) \ll N^{3/4+\varepsilon} (|\alpha| + |\beta| + |\gamma|)^{3/4 + \varepsilon},
% \end{equation}
% the implied constant depending on $\varepsilon$ only.
% \end{mytheo}

We also require the following result 
\begin{mytheo}
\label{theo:Molteni}
Let notation be as in this section and suppose $F$ is a Hecke-Maass form for $SL_3(\mz)$ (not necessarily arising as a symmetric-square lift).
Then for any $\varepsilon > 0$ we have
\begin{equation}
\label{eq:RankinSelbergsumA}
 \sum_{n \leq x} |A_F(1,n)|^2 \ll_{\varepsilon} x^{1 + \varepsilon} \lambda_F(\Delta)^{\varepsilon}.
\end{equation}
The implied constant is independent of $F$.  In the special case that $F$ is a symmetric-square lift of a $SL_2(\mz)$ Maass form, then 
% \begin{equation}
% \label{eq:RankinSelbergsumB}
%  \sum_{n \leq x} |A_F(1,n)|^3 \ll_{\varepsilon} x^{1+\varepsilon} T^{\varepsilon},
% \end{equation}
% and
\begin{equation}
 \label{eq:RankinSelbergsumC}
 \sum_{n \leq x} |A_F(1,n)|^4 \ll_{\varepsilon} x^{1+\varepsilon} T^{\varepsilon}
\end{equation}
\end{mytheo}
The crucial point in Theorem \ref{theo:Molteni} is the uniformity in terms of $F$.  Similar results to \eqref{eq:RankinSelbergsumA} and \eqref{eq:RankinSelbergsumC} without any explicit dependency on $F$, were given in \cite{LL}.  We give the proof of Theorem \ref{theo:Molteni} in Section \ref{section:Molteni}.  The proof is based on the following Lemma.
\begin{mylemma}
\label{lemma:8thmoment}
Let $\lambda_j(n)$ be the $n$-th Hecke eigenvalue of a Hecke-Maass form $u_j$ for $SL_2(\mz)$ with Laplace eigenvalue $1/4 + t_j^2$.  Then we have
\begin{equation}
\label{eq:8th}
\sum_{n \leq x} |\lambda_j(n^2)|^4 \ll x^{1+\varepsilon} t_j^{\varepsilon} \quad \text{and} \quad \sum_{n \leq x} |\lambda_j(n)|^8 \ll x^{1+\varepsilon} t_j^{\varepsilon},
\end{equation}
with implied constants depending on $\varepsilon > 0$ only.
\end{mylemma}

\begin{mycoro}
\label{coro:Moltenishort}
Suppose $B \leq A$ and $A > 0$.  Then uniformly in $F$,
\begin{equation}
 \sum_{n = A + O(B)} \frac{|A_F(1,n)|}{n} \ll (B/A)^{p} A^{\varepsilon} \lambda_F(\Delta)^{\varepsilon}.
\end{equation}
where we may take $p=1/2$ if $F$ is a general $SL_3(\mz)$ Maass form, $p=3/4$ if $F$ is a symmetric-square lift, and the Ramanujan conjecture allows $p=1$.
\end{mycoro}
This follows from Theorem \ref{theo:Molteni} using H\"{o}lder's inequality.

% Suppose that our progress towards the Ramanujan conjecture for a $GL_2$ Maass form leads to $|\lambda_j(p^2)| \leq 3 p^{\mu}$; so Kim-Sarnak (make ref) gives $\mu = 7/32$.  
% \begin{mylemma}
% \label{lemma:Ramanujanshort}
% For $B \leq A$ and $A > 0$ we have
% \begin{equation}
%  \sum_{n = A + O(B)} \frac{|A_F(1,n)|}{n} \ll \frac{B}{A} A^{\mu + \varepsilon}.
% \end{equation}
% \end{mylemma}

\section{Initial steps}
Let $Q \geq 1$ be a parameter to be chosen later.  By the Dirichlet approximation theorem, there exist coprime integers $a,q$ with $1 \leq q \leq Q$ so that $\alpha = \frac{a}{q} + \frac{\theta}{2\pi}$ with $|\frac{\theta}{2\pi}| \leq (qQ)^{-1}$.  Then with this notation the left hand side of \eqref{eq:smoothsum} takes the form
\begin{equation}
\label{eq:Sdef}
 S = \sum_{n \geq 1} A_F(1,n) e\big(\frac{an}{q} \big) \psi(n),
\end{equation}
with
\begin{equation}
\label{eq:psi}
\psi(y) = e^{i \theta y} w(y).
\end{equation}
Then with $d = \overline{a}$, (and $c$ switched with $q$) $S$ is precisely the left hand side of \eqref{eq:Voronoi}.  In order to analyze $S$ we then need to obtain satisfactory information on the behavior of the integral transforms $\Psi_{\pm}(x)$.  
It is convenient to record the following easy result.
\begin{mylemma}
\label{lemma:UsingWeil}
 Suppose that $\psi$ is given by \eqref{eq:psi} and $\Psi_{\pm}$ is defined by \eqref{eq:psidef}-\eqref{eq:Psi-}.  Then
\begin{equation}
\label{eq:PostVoronoi}
 |S| \ll q^{3/2+\varepsilon}\max_{\pm} \max_{d | q} \max_{n_1 | q/d} \sum_{n \geq 1} \frac{|A_F(n,1)|}{n} |\Psi_{\pm}(\frac{n n_1^2}{(q/d)^3})|,
\end{equation}
uniformly in $F$.
\end{mylemma}
This is a simple application of Weil's bound for Kloosterman sums combined with the Hecke relations.  In practice %, our bounds on $\Psi_{\pm}(x)$ shall only improve with larger values of $n_1$ and $d$ so that 
the important case is $d=n_1=1$.
\begin{proof}
By the Voronoi formula and Weil's bound, the left hand side of \eqref{eq:PostVoronoi} is 
\begin{equation}
\label{eq:Sbound}
\ll \max_{\pm} q \sum_{n_1 | q} \sum_{n_2 > 0} \frac{|A_F(n_2, n_1)|}{n_2 n_1} (q/n_1)^{1/2} d(q) |\Psi_{\pm}\leg{n_2 n_1^2}{q^3}|,
\end{equation}
where $d(q)$ is the divisor function.  
Next we use the Hecke relation $A_F(n_2, n_1) = \sum_{d |(n_2, n_1)} \mu(d) A_F(n_2/d, 1) A_F(1, n_1/d)$ (see Theorem 6.4.11 of \cite{Goldfeld}) followed by the triangle inequality and the bound $|A_F(1, l)| \ll l^{1/2}$ (uniformly in $F$), we see that $|A_F(n_2, n_1)| \ll \sum_{d | (n_1, n_2)} |A_F(n_2/d, 1)| (n_1/d)^{1/2}$.
Using this bound in \eqref{eq:Sbound}, and reversing the orders of summation, we have
\begin{equation}
 |S| \ll \max_{\pm} q^{3/2 + \varepsilon} \sum_{d|q} \sum_{n_1 | q/d} \frac{1}{n_1 \sqrt{d}} \sum_{n_2 > 0} \frac{|A_F(n_2, 1)|}{n_2}  |\Psi_{\pm}\leg{n_2 n_1^2}{(q/d)^3}|.
\end{equation}
Making $n = n_1 n_2$ be a new variable and taking the max over $d$ and $n_1$ gives the desired bound.
\end{proof}

% To use Lemma \ref{lemma:UsingWeil}, all that remains is to understand the size of $\Psi_{\pm}(x)$, which takes up the rest of the paper.  Once this is done (see Lemma ), it is straightforward to show the following 
% two different bounds on $S$ depending on the sizes of $q$ and $\theta$.  We do this in Section ??
% \begin{mylemma}
% \label{lemma:smallbeta}
%  Suppose $|\theta N| \leq 1$.  Then
% \begin{equation}
%  S \ll q^{3/2} T (qTN)^{\varepsilon}.
% \end{equation}
% \end{mylemma}
% \begin{mylemma}
% \label{lemma:largebeta}
%  Suppose $|\theta N| > 1$.  If the Ramanujan conjecture holds for $F$, then
% \begin{equation}
%  S \ll [q^{3/2} T + (N|\theta| q)^{3/2} ] (QTN)^{\varepsilon}
% \end{equation}
% Unconditionally, we have
% \begin{equation}
%  ?? S \ll [q^{3/2} T^{3/2} + (N|\theta| q)^{3/2} ] (QTN)^{\varepsilon}
% \end{equation}
% \end{mylemma}
% Since $|\theta| q \ll Q^{-1}$, and $q \leq Q$, we can claim for all $\lambda$ assuming Ramanujan that 
% \begin{equation}
%  S \ll (N^{3/2} Q^{-3/2} + Q^{3/2}T) (QTN)^{\varepsilon}.
% \end{equation}
% Choosing $Q = N^{1/2} T^{-1/3}$ optimizes the bound (here we may assume $Q \geq 1$ since $T \ll N^{1/4}$ is required for Theorem \ref{thm:Millertype} to be nontrivial.  %Unconditionally, we obtain
% \begin{equation}
%  S \ll (N^{3/2} Q^{-3/2} + Q^{3/2}T^{3/2} ) (QTN)^{\varepsilon},
% \end{equation}
% and choosing $Q = N^{1/2} T^{-1/2}$ optimizes the bound (again $Q \geq 1$ for the same reason as above) and gives Theorem \ref{thm:Millertype2}, in the smoothed form.  
\section{A Fourier-Mellin transform}
As a first step in understanding $\Psi_{\pm}$, we require information on the behavior of $\widetilde{\psi}(s)$, with $\psi(x)$ given by \eqref{eq:psi}.  To this end, we have
\begin{mylemma}
\label{lemma:integraltransform}
 Suppose that $\tau$, $\theta$, and $N$ are real numbers, and $w$ satisfies \eqref{eq:wbounds}.
Let
\begin{equation}
 I = \int_0^{\infty} w(x) e^{i \theta x} x^{i \tau} \frac{dx}{x}.
\end{equation}
Suppose that $|\tau| \geq 1$ and $|\theta N| \geq 1$.  Then
%Then there exists a sequence of functions $w_k$ each satisfying \eqref{eq:wbounds} (with a different set of constants $c_j$ for each $k$), such that
\begin{equation}
\label{eq:Iasymp}
 I = \sqrt{2 \pi} w(-\tau/\theta) |\tau|^{-\half} e^{i\tau \log|\tau /(e\theta)|} e^{\frac{\pi i}{4} \text{sgn}(\theta)} + O(|\tau|^{-3/2}).
\end{equation}
Furthermore, if $|\tau| \geq |\theta N|^{1+\varepsilon}$ then
\begin{equation}
\label{eq:taubig}
 I \ll_{A, \varepsilon} |\tau|^{-A}.
\end{equation}
and if $|\tau| \leq |\theta N|^{1-\varepsilon}$ then
\begin{equation}
\label{eq:tausmall}
 I \ll_{A, \varepsilon} |\theta N|^{-A}.
\end{equation}

%Here $w_0 = w$, and if $\tau$ and $\theta$ have the same sign then the main terms vanish.
\end{mylemma}
{\bf Remark.}
The conditions that $|\tau| \geq 1$ and $|\theta N| \geq 1$ ensure that the integral is oscillatory.  If these inequalities do not hold then the behavior of $I$ is easily determined by the following reasoning.
If $|\tau| \leq 1$ then the function $w_{\tau}(x) := w(x) (x/N)^{-1+i \tau}$ satisfies \eqref{eq:wbounds} with different absolute constants, in which case
\begin{equation}
 I =  N^{-1+i\tau}\intR w_{\tau}(x) e^{i \theta x} dx = N^{-1+i\tau} \widehat{w_{\tau}}(-\theta/(2\pi)),
\end{equation}
which displays the essential behavior of $I$ in this case (in particular, $I \ll_A  |\theta N|^{-A}$).  Similarly, if $|\theta N| \leq 1$ then $w_{\theta}(x) := w(x) e^{i \theta x}$ satisfies \eqref{eq:wbounds} with different constants, and the behavior of $I$ is fully determined by
\begin{equation}
 I = \widetilde{w_{\theta}}(i \tau),
\end{equation}
where $\widetilde{w_{\theta}}$ is the Mellin transform of $w_{\theta}$.  
In particular,
\begin{equation}
\label{eq:Ibound}
 I \ll_{A} (1 + |\tau|)^{-A}.
\end{equation}
\begin{proof}
First we show \eqref{eq:taubig} assuming $|\tau| \geq |\theta N|^{1+\varepsilon}$.  Let $g(x) = w(x) e^{i \theta x} $ so that $I = \widetilde{g}(i\tau)$.  Integrating by parts $j$ times shows $I \ll (|\theta N|/|\tau|)^j$; choosing $j$ large compared to $\varepsilon$ gives the desired bound.  The bound \eqref{eq:tausmall} is very similar but with $g(x) = w(x) x^{-1+i\tau}$, so that $I = \widehat{g}(-\theta/2\pi)$.  Integrating by parts $j$ times shows $I \ll (|\tau|/|\theta N|)^{j}$ which quickly gives \eqref{eq:tausmall} taking $j$ large.

Now we show \eqref{eq:Iasymp}.  We may assume $|\theta N|^{1-\varepsilon} \leq |\tau| \leq |\theta N|^{1+\varepsilon}$ since otherwise \eqref{eq:Iasymp} is consistent with the above analysis.  Write $I = \int_0^{\infty} g(x) e(f(x)) dx$ where $g(x) = w(x) x^{-1}$ and $2 \pi f(x) = \theta x + \tau \log{x}$.  Then $2 \pi f'(x) = \theta + \frac{\tau}{x}$ and for $j \geq 2$ $|f^{(j)}(x)| \asymp \frac{|\tau|}{x^{j}} \asymp \frac{|\tau|}{N^j}$.  Note that the stationary point of $f$ is $x_0 = -\tau/\theta$.  Suppose first that $x_0$ is not near the support of $w$ (which recall is a subset of $[N, 2N]$) in the sense that $x_0 \leq N/2$ or $x_0 \geq 4N$.  In this case, $|f'(x)| \gg \max(|\theta|, |\tau/N|)$.  Furthermore, $g^{(s)}(x) \ll \frac{1}{N^{s+1}}$ for $s=0, 1, 2, \dots$.  Then we apply Lemma 5.5.5 of \cite{Huxley} with Huxley's notation $(T, M, U, N, \alpha, \beta)$ taking the values in our notation $(|\tau|, 100N, N^{-1}, N/2, 2N)$.  Then we easily read off the bound $I \ll N^{\varepsilon} |\tau|^{-2}$ which is consistent with \eqref{eq:Iasymp} (note that the ``main term'' vanishes in this case).

Next suppose that $N/2 \leq x_0 \leq 4N$, in which case it is appropriate to 
apply Lemma 5.5.6 of \cite{Huxley} which gives the asymptotic formula for a weighted stationary phase integral.  In this case we have Huxley's parameters $(T, M, U, N, \alpha, \beta)$ taking the values in our notation $(|\tau|, 100N, N^{-1}, N, N/8, 8N)$.  We chose $\alpha, \beta$ so that $(\beta - x_0) \asymp (x_0 -\alpha) \asymp N$ so that the stationary point is not close to the endpoint of the range of integration.  The error term in Lemma 5.5.6 of \cite{Huxley} is then calculated to be $O(|\tau|^{-3/2})$.  If $\tau < 0$ then $f'(x)$ changes sign from negative to positive at $x_0 = -\tau/\theta$ and $f''(x) > 0$ so that we may directly compute the main term in Huxley's lemma, giving the stated main term in \eqref{eq:Iasymp}.  If $\tau > 0$ then Huxley's lemma applies directly to the complex conjugate of $I$ which after some easy manipulations leads to \eqref{eq:Iasymp} in this case too.
\end{proof}

\section{Bounds on $\Psi_{\pm}(x)$}
\label{section:smallbeta}
\begin{mylemma}
\label{lemma:phibound1}
Let $\Psi = \Psi_{k}$ (with either choice of $k=0$ or $1$) defined by \eqref{eq:psikdef}, with $\psi(x)$ defined by \eqref{eq:psi}.  Let
\begin{equation}
\label{eq:U}
 U = \max\{1+ T^2, T^2 |\theta N|, |\theta N|^3\},
\end{equation}
and
\begin{equation}
\label{eq:Delta}
\Delta = |xN - \frac{1}{(2 \pi)^3} |\theta N| T^2|.
\end{equation}
Then
% If $|\theta N| \leq T^{\varepsilon}$ then for any $A > 0$ we have
% \begin{equation}
% \label{eq:Phibound1}
%  \Psi(x) \ll_A T^{1+\varepsilon} (1 + \frac{xN}{T^{2+\varepsilon}})^{-A}
% \end{equation}
% with an implied constant depending on $w$ and $A$ only.  If $|\theta N| \geq T^{\varepsilon}$ then with 
\begin{equation}
\label{eq:Psibound}
 \Psi(x) \ll \mathcal{M} + \mathcal{E}_{\Delta},
\end{equation}
where
\begin{equation}
\label{eq:M}
 \mathcal{M} = %(T^2 + |\theta N|^2)^{1/2 + \varepsilon} |\theta N|^{1/2 + \varepsilon} 
\max(1+T, |\theta N|^{3/2}) | N T|^{\varepsilon}
(1 + \frac{xN}{U (NT)^{\varepsilon}})^{-A}
\end{equation}
and with certain absolute implied constants,
\begin{equation}
\label{eq:E}
 \mathcal{E}_{\Delta} = 
\begin{cases}
\frac{T^2}{|\theta N|^{1/2}}, \qquad &\text{ if }  T^{2/3} \leq |\theta N| \leq T^{1-\varepsilon} \text{ and } \Delta \ll |\theta N|^3, \\
\frac{T^3 |\theta N|}{\Delta}, \qquad &\text{ if } T^{2/3} \leq |\theta N| \leq T^{1-\varepsilon} \text{ and } |\theta N|^3 \ll \Delta \ll |\theta N| T^2, \\
%|\theta N | T, &\text{ if } T^{\varepsilon} \leq |\theta N| \leq T^{2/3} \text{ and } \Delta \ll T^2, \\ 
|\theta N| T \min(1, \frac{T^2}{\Delta}), &\text{ if } T^{\varepsilon} \leq |\theta N| \leq T^{2/3} \text{ and } \Delta  \ll |\theta N| T^2, \\
0, \qquad & \text{otherwise}.
\end{cases}
\end{equation}
% \begin{equation}
%  \label{eq:phibound3}
% \Psi(x) \ll |\theta N|^{3/2 + \varepsilon} (1 + \frac{xN}{|\theta N|^{3+\varepsilon}})^{-A}.
% \end{equation}
% If $T^{\varepsilon} \leq |\theta N| \leq T^{2/3}$ then
% \begin{equation}
%  \Psi(x) \ll |\theta N| T \delta_x + ??,
% \end{equation}
% where $\delta_x$ is the characteristic function of the interval $|xN - ?|\theta N| T^2| \ll T^2$.
% Finally, if $T^{2/3} \leq |\theta N| \leq T^{1-\varepsilon}$ then
% \begin{equation}
% \label{eq:phibound4}
%  \Psi(x) \ll \frac{T^2}{\sqrt{|\theta N|}} \delta'_x + |\theta N|^{3/2 + \varepsilon} (1 + \frac{xN}{|\theta N|^{3+\varepsilon}})^{-A} + \frac{T^3 |\theta N|}{?},
% \end{equation}
% where $\delta_x'$ is the characteristic function of the interval $|xN - ?|\theta N| T^2| \ll |\theta N|^3$.
\end{mylemma}

A heuristic calculation shows that if $\Delta \asymp |\theta N|^3 \gg T^2$ then the bound expressed by \eqref{eq:E} is essentially sharp since it is consistent with the asymptotic arising from stationary phase.  Likewise, if $\Delta \asymp T^{2-\varepsilon}$ and $|\theta N| \leq T^{2/3-\varepsilon}$ the bound is also essentially sharp and there is no cancellation in the integral.  In other cases there could potentially be extra savings over what is stated in \eqref{eq:E} by repeatedly integrating by parts but we did not investigate this.

Thus $x$ does indeed become localized in a relatively short interval (at least for certain ranges of the parameters).
It is therefore unavoidable that we study sums of the $GL_3$ Fourier coefficients in short intervals.  This is why the Ramanujan conjecture enters the picture.  Also note that if the $GL_3$ Maass form is fixed (as in \cite{Miller}) then this short-interval behavior is not really present.
% Lemma \ref{lemma:phibound1} combined with Lemma \ref{lemma:UsingWeil} immediately gives Lemma \ref{lemma:smallbeta}.
\begin{proof}
Recall that
\begin{equation}
 \widetilde{\psi}(-\sigma + i\tau) = \int_0^{\infty} w(x) x^{-\sigma} e^{i \theta x} x^{i \tau} \frac{dx}{x}.
\end{equation}
If $|\theta N| \leq 1$, then a small variation on \eqref{eq:Ibound} then shows
\begin{equation}
\widetilde{\psi}(-\sigma + i\tau) \ll_{A, \sigma} N^{-\sigma} (1+ |\tau|)^{-A},
% \frac{d^j}{d \tau^j} \widetilde{\phi}(-\sigma + i\tau) \ll_{A,j} N^{-\sigma} (\log{N})^j (1+ |\tau|)^{-A},
\end{equation}
with the implied constant independent of $|\theta|$ and $N$.  If $|\theta N| >1$ then Lemma \ref{lemma:integraltransform} applies.  We unify the cases with the bound $\widetilde{\psi}(-\sigma + i\tau) \ll_{\sigma, \varepsilon, A} N^{-\sigma} (1 + \frac{|\tau|}{1 + |\theta N|^{1+\varepsilon}})^{-A}$, valid for all $\theta$ and any $A > 0$.
By Stirling's approximation, and for $k=0$ or $1$,
\begin{multline}
 \left|
\frac{\Gamma\leg{1+\sigma -i\tau + iT+k}{2}}{\Gamma\leg{-\sigma + i \tau-iT+k}{2}} 
\frac{\Gamma\leg{1+\sigma -i\tau +k}{2}}{\Gamma\leg{-\sigma + i \tau+k}{2}} 
\frac{\Gamma\leg{1+\sigma -i\tau -iT+k}{2}}{\Gamma\leg{-\sigma + i \tau+iT+k}{2}}  
\right|
 \ll_{\sigma} [(1+|\tau - T|)(1+ |\tau|)(1+|\tau+T|)]^{\sigma + \half}
\\
\ll (|\tau|^3 + (|\tau| + 1) T^2)^{\sigma + \half}.
\end{multline}
Thus we obtain for $\sigma > -1$, 
\begin{multline}
\label{eq:PsiIntegralBound}
 \Psi(x) \ll_{\sigma, A} \intR (xN)^{-\sigma} (1 + \frac{|\tau|}{1 + |\theta N|^{1+\varepsilon}}))^{-A} (|\tau|^3 + (|\tau|+1) T^2)^{\sigma + \half} d\tau 
\\
\ll (1+|\theta N|^{1+\varepsilon}) U^{1/2} \leg{U}{xN}^{\sigma}.
\end{multline}
In the special case $|\theta N| \leq T^{\varepsilon}$, then \eqref{eq:PsiIntegralBound} is satisfactory for \eqref{eq:M} by choosing $\sigma = 0$ if $xN \leq U (NT)^{\varepsilon}$, and $\sigma$ large otherwise.  In this case the integration over $\tau$ is very short so that this is essentially optimal (except possibly for small values of $x$ which are not important in our application).

In fact, if $xN \geq U (NT)^{\varepsilon}$ then \eqref{eq:PsiIntegralBound} is satisfactory for \eqref{eq:M} for any range of $|\theta N|$.  
For the rest of the proof we therefore assume $|\theta N| \gg T^{\varepsilon}$, 
\begin{equation}
\label{eq:xNbound}
xN \leq U (NT)^{\varepsilon}, 
\end{equation}
and we fix $\sigma = -\half$ (for convenience).  Since $\widetilde{\psi}(-\sigma + i \tau)$ is very small except for $|\theta N|^{1-\varepsilon} \ll |\tau| \ll |\theta N|^{1+\varepsilon}$, we may restrict $\tau$ to such an interval, in the definition of $\Psi$; the error so obtained is much smaller than what is to be shown.  In this region of $\tau$, we use the asymptotic formula \eqref{eq:Iasymp} in which the main term takes the form
\begin{equation}
 N^{\half} W\leg{-\tau}{\theta} |\tau|^{-\half} e^{i\tau \log|\tau /(e\theta)|},
\end{equation}
where $W$ is a function satisfying \eqref{eq:wbounds}.  In particular, the support on $W$ implies $\tau$ has the opposite sign of $\theta$, and $|\tau| \asymp |\theta N|$.   For $k=0,1$, let $\Phi_k$ be given by the integral formula \eqref{eq:psikdef} but with $\widetilde{\psi}(-s)$ replaced by the main term from \eqref{eq:Iasymp}.  Precisely, let
\begin{multline}
 \Phi_k(x) = -\frac{(xN \pi^{3})^{1/2}}{2 \pi} \intR (x \pi^{3})^{i \tau} W\leg{-\tau}{\theta} |\tau|^{-\half} e^{i\tau \log|\tau /(e\theta)|}
\\
\frac{\Gamma\leg{1+\sigma -i\tau-iT + k}{2}}{\Gamma\leg{-\sigma + i \tau-iT+k}{2}} \frac{\Gamma\leg{1+\sigma -i\tau +k}{2}}{\Gamma\leg{-\sigma + i \tau+k}{2}}  \frac{\Gamma\leg{1+\sigma -i\tau+iT +k}{2}}{\Gamma\leg{-\sigma + i \tau+iT+k}{2}} d\tau.
\end{multline}
The error term satisfies
\begin{equation}
|\Psi_{k}(x) - \Phi_{k}(x)| \ll \frac{\sqrt{xN}}{(|\theta N| + T)^{100}} + \sqrt{xN} \int_{|\theta N|^{1-\varepsilon} \ll |\tau| \ll |\theta N|^{1+\varepsilon}} |\tau|^{-3/2}  d\tau \ll \frac{\sqrt{xN}}{|\theta N|^{1/2-\varepsilon}},
\end{equation}
which is consistent with \eqref{eq:M}.

We need to understand the oscillatory behavior of the gamma factors to determine the behavior of $\Phi_k$.  Stirling's approximation gives for $|t| \rightarrow \infty$
\begin{equation}
 \frac{\Gamma\leg{1+\sigma -it+k}{2}}{\Gamma\leg{-\sigma +it+k}{2}} =  |t/2|^{\sigma + \half} e^{- i t \log|t/2e|} (c_0 + \frac{c_1}{|t|} + \dots + O(\frac{1}{|t|^A}) ),
\end{equation}
where the $c_i$ are constants depending only on $k$ and the sign of $t$.
Write $\Phi_k(x) = \Phi_{k1}(x) + \Phi_{k2}(x)$, where $\Phi_{k2}$ corresponds to the portions of the integral with either $|\tau - T| \leq \sqrt{T}$ or $|\tau + T| \leq \sqrt{T}$.  A trivial bound gives
\begin{equation}
\label{eq:Phi2bound}
 \Phi_{k2}(x) \ll \sqrt{xN},
\end{equation}
which is satisfactory for \eqref{eq:M} upon noting that $|\theta N| \asymp T$ in this situation.

Then using Stirling's approximation, we obtain an asymptotic expansion for $\Phi_{k1}(x)$ as a sum of expressions of the form $\sqrt{xN} J$, where
\begin{equation}
\label{eq:Jdef}
J= \int_{|\tau \pm T| > \sqrt{T}}  g(\tau) e^{i f(\tau)} d\tau,
\end{equation}
where 
\begin{equation}
 f(\tau) = \tau \log(\frac{x |\tau| \pi^{3}}{e| \theta| } )  - (T+\tau)\log(|T+\tau|/2e) - \tau \log|\tau/2e| - (\tau-T) \log(|\tau-T|/2e),
\end{equation}
and $g(\tau)$ is a smooth function with support in an interval with $|\tau| \asymp |\theta N|$, and satisfying bounds of the form
\begin{equation}
 \frac{d^j}{d \tau^j} g(\tau) \ll |\tau|^{-\half -j}.
\end{equation}
Furthermore, the error in this expansion can be made to be $O(T^{-A})$ for $A$ arbitrarily large.  It therefore suffices to obtain upper bounds on $J$.

Notice that $f(\tau)$ simplifies a bit as
\begin{equation}
 f(\tau) = \tau \log(\frac{2 \pi^{3} x}{| \theta| } )  - (T+\tau)\log(|T+\tau|/2e)  - (\tau-T) \log(|\tau-T|/2e).
\end{equation}
We compute the derivatives of $f$:
\begin{equation}
\label{eq:f'}
 f'(\tau) = \log\Big(\frac{8\pi^3 xN}{ |\theta N(T+\tau)(T-\tau)|} \Big),
\end{equation}
and
\begin{equation}
 f''(\tau) = -\frac{1}{T+\tau} - \frac{1}{\tau-T} = - \frac{2 \tau}{\tau^2-T^2}.
\end{equation}

Now we treat cases.  Suppose $|\theta N| \geq T^{1-\varepsilon}$, so $U \ll  |\theta N|^3 T^{\varepsilon}$.  We recall from say Lemma 5.1.3 \cite{Huxley} that $\int_{\alpha}^{\beta} g(\tau) e^{if(\tau)} d\tau \ll V/\sqrt{\lambda}$ where $V$ is the total variation of $g$ along the interval of integration, plus the maximum modulus of $g$, and where $|f''(\tau)| \gg \lambda > 0$ along the interval.  In our case, $V \asymp |\theta N|^{-\half}$, and $\lambda \gg |\theta N|^{-1-\varepsilon}$, showing the desired bound $J \ll T^{\varepsilon}$ here.

Now suppose that $T^{2/3} \leq |\theta N| \leq T^{1-\varepsilon}$, so $U = |\theta N| T^2$.  We recall from say Lemma 5.1.2 of \cite{Huxley} that $\int_{\alpha}^{\beta} g(\tau) e^{if(\tau)} \ll V/\kappa$ where $V$ is the total variation of $g$ along the interval of integration, plus the maximum modulus of $g$, and where $\kappa$ is the infimum of $|f'|$ along the interval.
Notice that $f'(\tau) = \log(\frac{8\pi^3 xN}{ |\theta N| T^2}) + \frac{\tau^2}{T^2}(1+o(1))$.  Thus unless $xN \asymp |\theta N| T^2$ then $|f'(\tau)| \gg 1$ and so the first derivative bound would show $J \ll |\theta N|^{-1/2}$ which is consistent with \eqref{eq:M}.  So suppose $xN \asymp |\theta N| T^2$.  If it is the case that
\begin{equation}
\label{eq:smalllog}
 |\log(\frac{8 \pi^3 xN}{ |\theta N| T^2})| \leq 100 \frac{|\theta N|^2}{T^2}
\end{equation}
then a Taylor expansion of the logarithm shows that $xN$ lies in an interval of the form $\frac{1}{(2 \pi)^3} |\theta N| T^2 + O(|\theta N|^3)$, that is, $\Delta \ll |\theta N|^3$.  For such values of $x$ we apply the van der Corput bound (again, Lemma 5.1.3 of \cite{Huxley}) giving $J \ll \frac{T}{\sqrt{|\theta N|}} \frac{1}{|\theta N|^{1/2}} \asymp \frac{T}{|\theta N|}$, consistent with \eqref{eq:E}.  In the other case where \eqref{eq:smalllog} does not hold then $|f'(\tau)| \asymp |\log(\frac{8xN}{\pi^3 |\theta N| T^2})| \asymp \frac{\Delta}{|\theta N| T^2}$ for all $\tau$ in the region of integration and so the first derivative bound (Lemma 5.1.2 of \cite{Huxley}) shows 
\begin{equation}
\label{eq:Jfirstderiv}
J \ll \frac{|\theta N| T^2}{\Delta} \frac{1}{|\theta N|^{1/2}},
\end{equation}
which is consistent with \eqref{eq:E}.

Finally, consider the range $T^{\varepsilon} \leq |\theta N| \leq T^{2/3}$.  If $\Delta \leq 100 T^2$ then one can observe that the integral \eqref{eq:Jdef} is not oscillatory so we claim only the trivial bound which gives \eqref{eq:E}.  If $\Delta > 100T^2$ then the bound \eqref{eq:Jfirstderiv} carries over to this case also.

We have treated all possible cases and shown bounds consistent with \eqref{eq:Psibound}, so the proof is complete.
\end{proof}

% Now we proceed to prove Lemma \ref{lemma:smallbeta}.  The basic idea is to use the same approach as in Section \ref{section:smallbeta}, but there is a difficulty with short intervals.  First suppose that $|\theta N| \gg T$ so that we get a contribution to $S$ of size
% \begin{equation}
% q^{3/2 + \varepsilon} |\theta N|^{3/2 + \varepsilon},
% \end{equation}
% consistent with Lemma \ref{lemma:smallbeta}.  
% 
% Suppose now that $|\theta N| \ll T$.
% Next suppose that $x = c \theta T^2 + O(\Delta)$ with $\Delta \gg \theta^3 N^2$.  In this range we get a contribution to $S$ of size
% \begin{equation}
% q^{3/2 +\varepsilon} T \frac{\theta T^2}{\Delta} \sum_{l = q^3 c \theta T^2 + O(q^3 \Delta)} \frac{|a_l|}{l},
% \end{equation}
% which upon using the Ramanujan bound gives a total bound of size
% \begin{equation}
% T q^{3/2 + \varepsilon} \ll |\theta N| q^{3/2 + \varepsilon}.  
% \end{equation}
% Note that the left hand bound above is much stronger than what we need so perhaps it's easy to remove Ramanujan here.
% If $\Delta \ll \theta^3 N^2$ then we get a bound to $S$ of size
% \begin{equation}
% T^{3/2} q^{3/2} \sum_{l = c\theta T^2 q^3 + O(q^3 \Delta)} \frac{|a_l|}{l},
% \end{equation}
% which again upon using Ramanujan is
% \begin{equation}
% \ll (qT)^{3/2} \frac{\Delta}{\theta T^2} \ll  (qT)^{3/2} \frac{\theta^2 N^2}{T^2} =(q\theta N)^{3/2} \leg{\theta N}{T}^{1/2} \ll (q\theta N)^{3/2},
% \end{equation}
% again consistent with Lemma \ref{lemma:largebeta}.

\section{Bounding $S$}
\subsection{A general bound}
Recall the definition of $S$ given by \eqref{eq:Sdef}.  Here we prove a bound on $S$ that is explicit in terms of $q$ and $\theta$.  In the next section we find a bound that is uniform in $\alpha$ and optimize the parameter $Q$.

Applying Lemma \ref{lemma:UsingWeil} and Lemma \ref{lemma:phibound1}, we obtain say $S \ll S_{\mathcal{M}} + S_{\mathcal{E}}$, corresponding to $\Psi \ll \mathcal{M} + \mathcal{E}$.
An easy calculation shows
\begin{equation}
 S_{\mathcal{M}} \ll q^{3/2} (T + |\theta N|^{3/2})(NTq)^{\varepsilon}.
\end{equation}
Using $q|\theta| \leq 2\pi Q^{-1}$, we have
\begin{equation}
\label{eq:SM}
 S_{\mathcal{M}} \ll (Q^{3/2} T + N^{3/2} Q^{-3/2})(NTQ)^{\varepsilon}.
\end{equation}

To bound $S_{\mathcal{E}}$, we break it into cases.  We show that if $T^{2/3} \leq |\theta N| \leq T^{1-\varepsilon}$ then
\begin{equation}
\label{eq:SE1}
 S_{\mathcal{E}} \ll q^{3/2} \frac{T^2}{|\theta N|^{1/2}} \leg{|\theta N|^2}{T^2}^p (NTq)^{\varepsilon}
%q^{3/2}  \min(T|\theta N|^{1/2}, |\theta N|^{3/2} (q^3 |\theta| T^2)^{\mu} (NTq)^{\varepsilon}
\end{equation}
where $p$ is as in Corollary \ref{coro:Moltenishort}, 
and if $T^{\varepsilon} \leq |\theta N| \leq T^{2/3}$ then
\begin{equation}
\label{eq:SE2}
 S_{\mathcal{E}} \ll q^{3/2} |\theta N|^{1-p} T (NTq)^{\varepsilon}.
%q^{3/2}  \min(T|\theta N|^{1/2}, T (q^3 |\theta| T^2)^{\mu} (NTq)^{\varepsilon}.
\end{equation}
Suppose first that $T^{2/3} \leq |\theta N| \leq T^{1-\varepsilon}$ and consider the contribution from $\Delta= |xN - \frac{1}{(2 \pi)^3} |\theta N| T^2| \ll |\theta N|^3$; recall $x = n n_1^2 d^3/q^3$.  We require Corollary \ref{coro:Moltenishort} with $A = \frac{q^3}{d^3 n_1^2} |\theta| T^2$ and $B= \frac{q^3}{d^3 n_1^2} |\theta|^3 N^2$.  
This gives a contribution to $S_{\mathcal{E}}$ of size
% \begin{equation}
% \label{eq:SE1case1}
% q^{3/2} \frac{T^2}{|\theta N|^{1/2}} \leg{|\theta N|^2}{T^2}^{3/4} (NTq)^{\varepsilon}
% \end{equation}
% which simplifies to 
given by \eqref{eq:SE1}.  It turns out that the case $|\theta N|^3 \ll \Delta \ll |\theta N| T^2$ is similar.  Say $Y \leq \Delta \leq 2Y$. Then we can dissect the interval $[Y, 2 Y] \subset \mr$ into $\ll Y |\theta N|^{-3}$ subintervals of length at most $|\theta N|^3$.  This leads to $O(Y |\theta N|^{-3})$ instances of bounds of the form \eqref{eq:SE1} but with $T^3 |\theta N|/Y$ in place of the term $T^2 |\theta N|^{-1/2}$ appearing in \eqref{eq:SE1}.  This leads to a bound of size
\begin{equation}
 q^{3/2} \frac{T^3}{|\theta N|^2}  \leg{|\theta N|^2}{T^2}^{p} (NTq)^{\varepsilon}.
\end{equation}
Since $T^3 |\theta N|^{-2} \leq T^2 |\theta N|^{-1/2}$ for this range on $|\theta N|$, the above bound is no worse than \eqref{eq:SE1}, as desired.  Thus we have proved \eqref{eq:SE1}.

The case of $T^{\varepsilon} \leq |\theta N| \leq T^{2/3}$ is similar.  We have in this case $A = \frac{q^3}{d^3 n_1^2} |\theta| T^2$ (as before), and $B = \frac{q^3}{d^3 n_1^2} \frac{T^2}{N}$, so $B/A = |\theta N|^{-1}$.  Then the contribution from $\Delta \ll T^2$ gives
\begin{equation}
 q^{3/2} |\theta N| T \frac{1}{|\theta N|^{p}} (NTq)^{\varepsilon},
\end{equation}
consistent with \eqref{eq:SE2}.  As in the previous case, the ranges with $T^2 \ll \Delta \ll |\theta N| T^2$ can be dissected into subintervals of length $T^2$ giving the same final contribution.  This proves \eqref{eq:SE2}.

\subsection{Uniform bound}
Our bound on $\mathcal{S}$ depending on $\theta$ and $q$ is given by adding \eqref{eq:SM} and possibly \eqref{eq:SE1} or \eqref{eq:SE2} if applicable.

First, suppose that the Ramanujan conjecture holds, so $p=1$.  In this case we obtain
\begin{equation}
\label{eq:Sp1}
 S \ll (Q^{3/2} T + N^{3/2} Q^{-3/2})(NTQ)^{\varepsilon}
\end{equation}
which leads to \eqref{eq:mainthmRama} after choosing $Q = N^{1/2} T^{-1/3}$ (observe that the bound in \eqref{eq:mainthmRama} is trivial if $N \leq T^2$ so that $Q \geq 1$).

Next we take the unconditional bound with $p=3/4$.  This leads to the same bound as $\eqref{eq:Sp1}$ coming from the contribution of $S_{\mathcal{M}}$, plus two extra terms, the one arising from \eqref{eq:SE1} having size
\begin{equation}
 q^{3/2} T^{1/2} |\theta N|  (NTq)^{\varepsilon} \ll (Q^{-1/2} T^{1/2} N)(NTQ)^{\varepsilon}, 
\end{equation}
and the other arising from \eqref{eq:SE2} of size
\begin{equation}
 Q^{3/2} T^{7/6} (NTQ)^{\varepsilon}.
\end{equation}
Taken together, these lead
to \eqref{eq:mainthmRama} after picking $Q = N^{1/2} T^{-1/3}$.

%********************************************************************************************
\section{Bounds on moments of Hecke eigenvalue}
\label{section:Molteni}
Here we give a proof of Theorem \ref{theo:Molteni}.
The bound \eqref{eq:RankinSelbergsumA} is a straightforward application of \cite{Brumley} while for \eqref{eq:RankinSelbergsumC} we require the convexity bound on $GL_4 \times GL_4$ due to \cite{XiannanLi}.  Most importantly, these bounds rely on the automorphy of the second through fourth symmetric powers \cite{Shimura} \cite{GJ} \cite{KimShahidi} \cite{Kim}.

Now we prove \eqref{eq:RankinSelbergsumC}.  This follows immediately from Lemma \ref{lemma:8thmoment}.
We shall show
\begin{equation}
\label{eq:S6bound}
S_{4,2}(x):=\sum_{n \leq x} |\lambda_j(n^2)|^4 \ll x^{1+\varepsilon} t_j^{\varepsilon},
\end{equation}
the other bound stated in \eqref{eq:8th} being similar.  Actually, we sketch an argument giving an asymptotic formula for the left hand side of \eqref{eq:S6bound} (say modified with a smooth weight) with a power saving in the error term with polynomial dependence on $t_j$; this is a modification of work in \cite{LL} where the dependency on $t_j$ is not explicated.

A simple Mellin transform argument using $e \cdot e^{-n/x} \geq 1$ for all $n \leq x$ then shows
\begin{equation}
 S_{4,2}(x) \leq e \frac{1}{2 \pi i} \int_{(2)} x^{s} \Gamma(s) D(s) ds,
\end{equation}
where $D(s) = \sum_{n=1}^{\infty} \lambda_j(n^2)^4 n^{-s}$.  We shall show that $D(s) \ll_{\delta, \varepsilon} t_j^{\varepsilon}$ for $s = \sigma + it$ and $\sigma \geq 1 +\delta$, $\delta > 0$ with the implied constant independent of the Maass form.  This immediately gives \eqref{eq:S6bound} after moving the contour of integration to $1 + \varepsilon$.  In fact we show that $D(s)$ has a meromorphic continuation to $\sigma > 23/32$ and satisfies a polynomial growth in terms of $t$ and $t_j$ in this region.

We shall express $D(s)$ in terms of $L$-functions associated to symmetric power lifts of $u_j$ as in \cite{LL}.  Writing $\lambda_j(p^k) = \sum_{m=0}^{k} \alpha_p^m \beta_p^{k-m}$ where $\alpha_p + \beta_p = \lambda_j(p)$ and $\alpha_p \beta_p = 1$, we obtain
\begin{equation}
 D(s) = \prod_{p} \sum_{k=0}^{\infty} (\sum_{m=0}^{2k} \alpha_p^m \beta_p^{2k-m})^4 p^{-ks}.
\end{equation}
Using Mathematica, say, we can compute this infinite series as a rational function in $p^{-s}$.  The main idea is to compare the coefficient of $p^{-s}$ with known $L$-functions (this can be done by hand).  In fact, we see that $D(s)$ has the same coefficient of $p^{-s}$ as 
\begin{equation}
E(s):=L(\text{sym}^4 u_j \times \text{sym}^4 u_j, s) L(\text{sym}^3 u_j \times \text{sym}^3 u_j, s)^2 L(\text{sym}^4 u_j, s)^3 L(\text{sym}^2 u_j, s)^3.
\end{equation}
Recall that
\begin{equation}
\label{eq:RSsym}
 L(\text{sym}^r u_j \times \text{sym}^r u_j, s) = \prod_{p} \prod_{m=0}^{r} \prod_{n=0}^{r} (1 - \alpha_p^{r-m} \beta_p^m \alpha_{p}^{r-n} \beta_p^{n} p^{-s})^{-1},
\end{equation}
and
\begin{equation}
 L(\text{sym}^r u_j, s) = \prod_{p} \prod_{m=0}^{r} (1 - \alpha_p^{r-m} \beta_p^m p^{-s})^{-1}.
\end{equation}
By work of \cite{Shimura} \cite{GJ} \cite{KimShahidi} \cite{Kim}, we have that $E(s)$ has a meromorphic continuation to $\mc$ with a pole at $s=1$ only.  Furthermore, Xiannan Li \cite{XiannanLi} showed that $E(s)$ satisfies the convexity bound, that is, $|E(1+ \delta + it)| \leq C(\delta, \varepsilon) t_j^{\varepsilon}$ for $\delta > 0$ and $C(\delta, \varepsilon)$ not depending on $t$ or $t_j$; inside the critical strip there is a polynomial bound in terms of $t_j$ and $s$.

Write $D(s) = E(s) U(s)$.  Then a computer calculation shows that $U(s) = \prod_p U_{p}(s)$ where
\begin{equation}
\label{eq:Ujpsbound}
 U_{p}(s) = 1 + O((|\alpha_p|^{12} + |\beta_p|^{12})p^{-2\sigma}).
\end{equation}
The implied constant is absolute.  Next observe that the convexity bound for $L(\text{sym}^4 u_j \times \text{sym}^4 u_j, s)$ implies that for $\sigma > 1 + \delta$, $\delta > 0$, we have
\begin{equation}
\label{eq:8thEP}
\prod_p (1 + \frac{|\alpha_p|^8 + |\beta_p|^8}{p^{\sigma}}) \leq C(\delta, \varepsilon) t_j^{\varepsilon}.
\end{equation}
The reason is that if $p$ is such that the Ramanujan bound holds, that is $|\lambda_j(p)| \leq 2$ then $|\alpha_p| = |\beta_p| = 1$
while if Ramanujan does not hold then $\alpha_p$ and $\beta_p$ are real and by positivity the Euler factor at $p$ in \eqref{eq:RSsym} with $r=4$ is larger than the above Euler factor at $p$.  Combining cases gives \eqref{eq:8thEP}.  We relate $U(s)$ to an instance of \eqref{eq:8thEP} by borrowing $\max(|\alpha_p|, |\beta_p|)^4$ and using the Kim-Sarnak bound \cite{Kim} $|\alpha_p|, |\beta_p| \leq p^{7/64}$ on this part.  Hence, we have
\begin{equation}
U(s) \ll \prod_p (1 + O( \frac{|\alpha_p|^{8} + |\beta_p|^{8}}{p^{2\sigma-\frac{7}{16}}} ) ).
\end{equation}
For $\sigma \geq \frac{23}{32} + \delta$ we then have a bound $U(s) \ll_{\delta, \varepsilon} t_j^{\varepsilon}$.

\section{Unsmoothing}
\label{section:unsmoothing}
Here we show how to deduce Theorem \ref{thm:Millertype} from Theorem \ref{thm:mainthmsmooth}.  The uniformity in $\alpha$ makes this deduction fairly easy with no essential losses in quality.  We begin with
\begin{mylemma}
\label{lemma:unsmoothing}
 For any integer $x \geq 1$ there exists a function $h(t)$ satisfying
\begin{equation}
\label{eq:hintbound}
 \intR |h(t)| dt \leq 7 + \log(x),
\end{equation}
and for any integer $n$,
\begin{equation}
 \intR h(t) e(nt) dt = \begin{cases}
                         1, \qquad |n| \leq x \\
			 0, \qquad \text{otherwise}.
                       \end{cases}
\end{equation}
\end{mylemma}
This is a simple variant on Lemma 9 of \cite{DFIbilinear} and we omit the proof.
% \begin{proof}
%  Let $f(v)$ be the even function defined for $v \geq 0$ by
% \begin{equation}
%  f(v) =
% \begin{cases}
%   1, \qquad &0 \leq v \leq x, \\
%   x+1-v, \qquad &x < v \leq x+1, \\
%   0, \qquad &v \geq x+1.
%   
% \end{cases}
% \end{equation}
% Then for integer $n$, $f(n) = 1$ if $|n| \leq x$ and is zero otherwise.  By the Fourier inversion theorem and the evenness of $f$, we have
% \begin{equation}
%  f(v) = \int_{-\infty}^{\infty} h(t) \cos(2\pi tv) dt, \qquad h(t) = 2 \int_0^{\infty} f(v) \cos(2 \pi vt) dv.
% \end{equation}
% A trivial estimation shows (a) $|h(t)| \leq 2(x+1)$.  Integration by parts shows
% \begin{align}
%  h(t) = \frac{2}{2 \pi t} \int_x^{x+1} \sin(2 \pi t v) dv = \frac{1}{2 \pi^2 t^2} (\cos(2 \pi tx) - \cos(2 \pi t(x+1)),
% \end{align}
% which give the further bounds (b) $|h(t)| \leq \frac{1}{\pi |t|}$ and (c) $|h(t)| \leq \frac{1}{\pi^2 |t|^2}$.  Using the bound (a) for $|t| \leq x^{-1}$, the bound (b) for $x^{-1} < |t| \leq 1$, and the bound (c) for $|t| > 1$ gives \eqref{eq:hintbound}.
% \end{proof}
Now we show how to prove Theorem \ref{thm:Millertype}.  Consider first the sum
\begin{equation}
 \sum_{2M/3 < n \leq M} A_F(1,n) e(\alpha n) = \sum_{2M/3 < n \leq M} A_F(1,n) e(\alpha n) w(n),
\end{equation}
where $w$ is a function satisfying \eqref{eq:wbounds}, with $M = 15N/8$ and $w(n) = 1$ for $2M/3 < n \leq M$.  Write the right hand side as $S_w(M)-S_w(2M/3)$.  Then apply Lemma \ref{lemma:unsmoothing} to get
\begin{equation}
 |S_w(x)| \leq \intR |h(t)| |\sum_n A_F(1,n) e(n(\alpha + t)) w(n) | dt .
\end{equation}
Theorem \ref{thm:mainthmsmooth} applies to the sum over $n$, the uniformity in $\alpha$ being critical, and \eqref{eq:hintbound} controls the $t$-integral, giving Theorem \ref{thm:Millertype}.

\end{document}